\RequirePackage{ifpdf}
\ifpdf
 	\documentclass[pdftex]{amsart}
 	\else
 	\documentclass[dvips]{amsart}
\fi

\usepackage{amsfonts,amsthm,latexsym,amsmath,amssymb,amscd,amsmath, mathrsfs, epsf}
\usepackage{graphicx}

\graphicspath{{./figures/}}
\ifpdf
	\usepackage{epstopdf}		
	\DeclareGraphicsExtensions{.png,.jpg,.eps,.epsf}
\fi

\newtheorem{theorem}{Theorem}

\newtheorem{lemma}[theorem]{Lemma}

\newtheorem{corollary}{Corollary}
\newtheorem{example}{Example}

\newtheorem{problem}{Problem}

 \newcommand{\eps}{\epsilon}
 \newcommand{\bC}{\mathbb{C}}

\newcommand{\Z}{\mathcal{Z}}

\newcommand{\al}{\alpha}

\begin{document}
\title[An asymptotic Gauss-Lucas theorem]{On asymptotic Gauss-Lucas theorem}

\author[R.~B\o gvad]{Rikard B\o gvad }

\address{ Department of Mathematics,
   Stockholm University,
   S-10691, Stockholm, Sweden}
\email{rikard@math.su.se}

\author[D.~Khavinson]{Dmitry Khavinson}
\address{Department of Mathematics, University of South Florida, Tampa,
         33620, USA}
\email{dkhavins@usf.edu}

\author[B. Shapiro]{Boris Shapiro}

\address{Department of Mathematics,
   Stockholm University,
   S-10691, Stockholm, Sweden}
\email{shapiro@math.su.se}

\begin{abstract}
In this note we extend the  Gauss-Lucas theorem on the zeros of the derivative of a univariate polynomial to the case of sequences of univariate polynomials  whose almost all  zeros  lie in a given convex bounded domain in $\bC$.
\end{abstract}

\maketitle

\section{Introduction}

The celebrated Gauss-Lucas theorem claims that for any univariate polynomial $P(z)$ with complex coefficients, all roots of $P^\prime(z)$ belong to the convex hull of the roots of $P(z)$, see  Theorem 6.1 of \cite {Ma}. Many generalizations have been obtained over the years see, e.g., \cite{Ch, CuMa, Ru} and references therein.

\medskip
In the present note, motivated by problems in potential theory in $\bC,$ we extend the Gauss-Lucas theorem to sequences of polynomials of increasing degrees whose almost all zeros  lie in a given convex bounded domain  in $\bC$.
 Namely, given a convex bounded domain $\Omega\subset \bC,$ let $\{p_n(z)\}_{n=0}^\infty$ be a sequence of univariate polynomials with the degrees
$\deg p_n=m_n$ such that $\lim_{n\to \infty} m_n=+\infty$. Assume that  $\lim_{n\to \infty}\frac{\sharp_n(\Omega)}{m_n}=1,$ where $\sharp_n(\Omega)$ is the  number of zeros of $p_n$  lying in $\Omega$ (counted with multiplicities).

\begin{problem}\label{prob:main}  Following the above notation we now ask whether there exists $\lim_{n\to \infty}\frac{\sharp_n^\prime(\Omega)}{m_n-1},$ where $\sharp_n^\prime(\Omega)$ denotes the number of zeros of $p_n^\prime(z)$ lying in $\Omega$?
\end{problem}

It turns out that the answer to  Problem~\ref{prob:main} formulated verbatim as above ,  is, in general, negative.

\begin{figure}

\begin{center}
\includegraphics[scale=0.4]{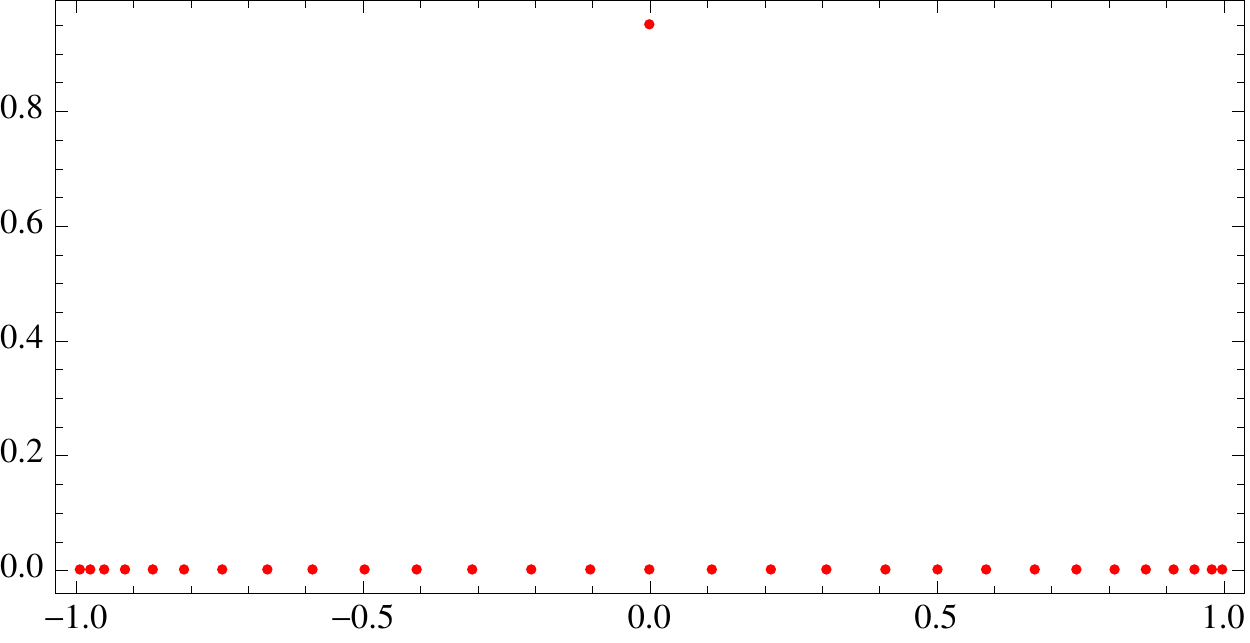}

\end{center}


\caption{Zeros of  $((z-i)T_{30}(z))'$}
\label{fig1}
\end{figure}

\begin{example}
{\rm Let $O$ be the open square $(-2,2)\times (-4i,0)$.
If $T_{n}(z):=\cos (n\arccos z)$ is the $n$-th Chebyshev polynomial of the first kind,
then the derivative of $(z-i)T_{n}(z)$ has all its
zeros in the upper half plane. Therefore,
if we replace $z$ by $z+ia_n$ for some sufficiently small $a_n$ (i.e.,  shift all zeros downward by $a_n$),
then we obtain polynomials of degrees $n+1$ with $n$ zeros in $O$,
but whose derivatives  have no zeros  in $O$. Choosing $a_n$ appropriately for each $n$, we get a sequence of polynomials with all but one zeros in $O$ whose derivatives have no zeros in $O$.

Strict convexity (e.g., as in the case of the open unit disk) will not be of much help either.
Just replace above $z$ by $M_nz$ with some large $M_n$ and then
make a vertical translation so that after all these operations
the image of $[-1,1]$ becomes a tiny secant
segment of the unit circle.  (This example was suggested to the third author by Professor V.~Totik.)}
\end{example}

\medskip
However  with  slightly weaker requirements  Problem 1 has a  positive answer.

\begin{theorem}\label{th:main} Given a polynomial sequence $\{p_n(z)\}$ as above and any $\eps>0$,
$$\lim_{n\to \infty}\frac{\sharp_n^\prime(\Omega_\eps)}{m_n-1}=1,$$ where $\sharp_n^\prime(\Omega_\eps)$ is the number of zeros of $p_n^\prime(z)$ lying in  $\Omega_\eps,$ and  $\Omega_\eps$ is the open $\eps$-neighborhood of $\Omega\subset \bC$.
\end{theorem}

Two illustrations of Theorem~\ref{th:main} are given in Figure~2.

\medskip
\noindent
{\it Acknowledgements.} The third author want to thank Professor V.~Totik for discussions of the problem.

\section{Proof}

\begin{figure}

\begin{center}
\includegraphics[scale=0.26]{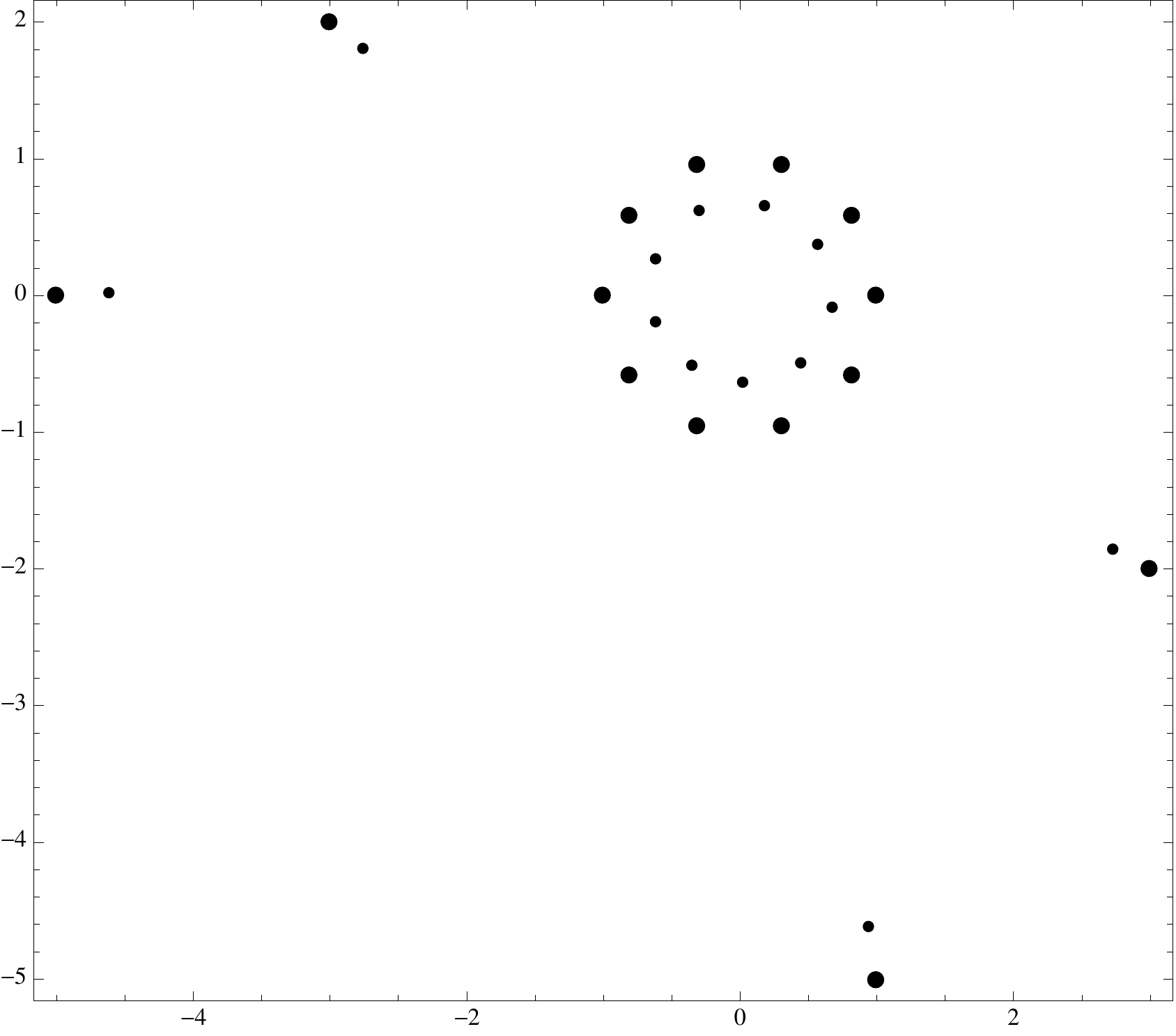} \hskip1cm \includegraphics[scale=0.2]{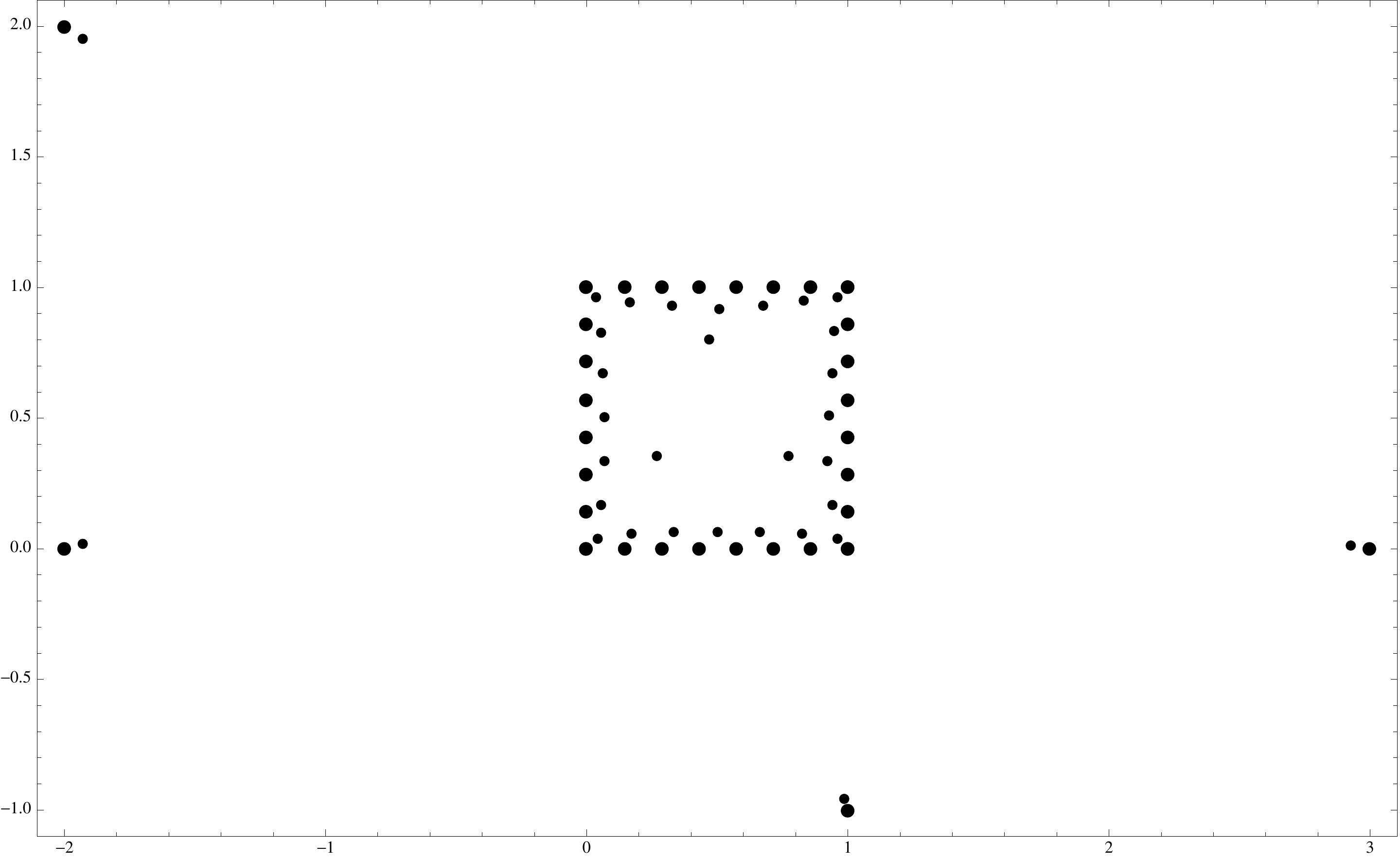}

\end{center}


\caption{Zeros of $P$ (larger dots) and $P'$ (smaller dots) for $P=(z^{10}-1)(z-3 + 2I)(z +3-2I)(z-1+5I)(z + 5)$ (left)
and  $P=Q(z-3)(z+2-2I)(z-1+I)(z+2)$ where $Q$ has $7$ uniformly placed zeros on each side of the unit square (right).}
\label{fig2}
\end{figure}

We first prove Theorem~\ref{th:main} in the case when $\Omega$ is a disk. 
Fix $\epsilon>0$.
 Let $\mathbb{D}=\{|z|<1\}$ be the open unit disk, $p_n(z)=\prod_{k=1}^{m_n}(z-\al_k)$, $\lim m_n=\infty$. Let us factor $p_n$ as follows:

\begin{equation}\label{(1)}
p_n=q_n r_n=\prod_{k=1}^{k_n}(z-a_k)\prod_{k=k_{n+1}}^{m_n}(z-b_k),\quad |b_k|>1+\epsilon, \quad |a_k|<1+\epsilon, \quad \text{with} \quad \lim_{n\to \infty} \frac{k_n}{m_n}=1.
\end{equation}

Denote by $\mathfrak{Z}_n':=\sharp\{z:|z|<1+\eps: p_n'(z)=0\}$. We want to show that $\lim_{n\to \infty} \frac{\mathfrak{Z_n'}}{m_n-1}=1$.

Let \begin{equation}\label{(2)}\hat \mu_n(z):=\frac{1}{m_n}\sum_{k=1}^{_n}\frac{1}{z-\al_k}=\frac{1}{m_n}\frac{p'_n}{p_n},
\end{equation}  denote the Cauchy transform of the root-counting probability measure $\mu_n$ of $p_n$. Note that

\begin{equation}\label{(3)}
\hat \mu_n(z)=\frac{1}{m_n}\left(\frac{q_n'}{q_n}+\frac{r_n'}{r_n}\right)=\frac{1}{m_n}(k_n\hat \nu_n+(m_n-k_n)\hat \psi_n),
\end{equation}
where $\nu_n$ and $\psi_n$ are the root-counting measures of $q_n$ and $r_n$ respectively. All zeros of $q_n'$ lie in the unit disk $D$ by the Gauss-Lucas theorem.  Also \eqref{(1)} implies that

\begin{equation}\label{(4)}
 \frac{(m_n-k_n)}{m_n}||\psi_n||\to 0.
\end{equation}

Formula \eqref{(4)} implies that for all $p: 1\leq p <2$, we have
\begin{equation}\label{(5)}
\left\vert\left \vert \frac{1}{m_n}\frac {r_n'}{r_m}\right\vert\right \vert_{L^p(dA)}\to 0
\end{equation}
 on compact subsets of $\bC$. Here, $dA=\frac{1}{\pi}dxdy$ denote the normalized area measure.

Equation \eqref{(5)} follows from a trivial observation. Let $\mu$ be a Borel measure with a compact support. Then for any compact set $K\subset \bC,\;$ and any $p: 1 \le p <2$, we have

\begin{equation}\label{(6)}
||\hat \mu(z)||^p_{L^p(K,dA)}\le C(p,K)||\mu||.
\end{equation}

Indeed, $\hat \mu=\int\frac{d\mu(\xi)}{\xi-z},$ hence
\begin{equation}\label{(7)}
\int_K|\hat \mu|^pdA\le ||\mu||\int_{|\xi|<R}\frac{1}{|\xi|^p}dA\le C||\mu||,
\end{equation}
where $R$ is chosen so that $\forall \xi\in K$ the disk of radius $R$ centered at $\xi$  contains $K$. The integral in \eqref{(7)} converges for all $p<2$ and \eqref{(6)} follows, hence does \eqref{(5)}.

Thus,  we have from \eqref{(5)} the following corollary.

\begin{corollary}\label{cor:1}
. For any fixed $R>1+\eps, \mbox{and \, any} \, p: 1\le p <2$ , for almost all $r:  1+\eps<r<R, $ we have
\begin{equation}\label{(8)}
\lim_{n\to \infty} \frac{1}{m_n}\int_{|z|=r} \frac{|r_n'|^p}{|r_n|^p}ds_r=0,
\end{equation}
where $ds_r$ is the arclength measure on $\{z: |z|=r\}$.
\end{corollary}

Thus, from \eqref{(3)}, we now obtain

\begin{corollary}
$$\lim_{n\to \infty} \frac{1}{m_n}\int_{|z|=r}\left\vert \frac{p_n'}{p_n}-\frac{q'_n}{q_n}\right\vert^pds_r=0$$
for almost all $r: 1+\eps<r<R$ and $p: 1\le  p<2.$
\end{corollary}

However
$$\frac{1}{m_n} \left( \frac{p_n'}{p_n}-\frac{q_n'}{q_n} \right)=\frac{1}{m_n} \left( \sum_{k=k_n+1}^{m_n} \frac{1} {z-b_k} \right)$$
by \eqref{(3)}, and hence is analytic inside $\{|z|<1+\epsilon\}$ since $|b_k|>1+\epsilon$.

Therefore, from standard results on Hardy spaces $H^p$ in the disk, cf. \cite{Du},we conclude that
\begin{equation}\label{(9)}
\frac {1} {m_n}  \left(  \frac{p_n'}{p_n}- \frac {q_n'} {q_n} \right) \to 0
\end{equation}
uniformly in the closed disk $\overline D=\{|z|\le 1+\epsilon$.

Since $\frac{1}{m_n}\frac{q_n'}{q_n}$ vanishes at $k_n-1$ points in $D$ by Gauss-Lucas theorem, invoking Hurwitz's theorem, we conclude that

\begin{equation}\label{(10)}\lim_{n\to\infty}\frac{1}{m_n}\left[ \sharp (z: |z|<1+\eps: p_n'(z)=0)-(k_n-1)\right]=0.
\end{equation}
Since, by assumption, $\lim_{n\to \infty} \frac {k_n-1} {m_n}=\lim_{n\to \infty}\frac{k_n}{m_n}=1$, we arrive at
 $$ \lim_{n\to \infty}\frac{\mathfrak{Z_n}'}{m_n}=\lim_{n\to \infty}\frac{\mathfrak{Z_n}'}{m_n-1}=1,$$
which settles Theorem~\ref{th:main} in the case of a disk.

To finish the proof for the general case of an arbitrary convex domain $\Omega$ observe that we only used some properties of a disk to get a convenient foliation of a neighborhood of the unit disk by concentric circles and when applying Gauss-Lucas theorem. Both these facts are readily available for an arbitrary bounded convex domain. Finally, the Hardy spaces of analytic functions in the disk are replaced by the Smirnov classes $E^p$ of functions representable by Cauchy integrals with $L^p$-densities (with respect to arclength). In the domains with piecewise smooth boundaries, e.g., convex domains, the latter behave in the very same manner as Hardy spaces -- cf. \cite{Du}.
\qed

\end{document}